%



\input amstex.tex

\magnification=\magstep1
\hsize=5.5truein
\vsize=9truein
\hoffset=0.5truein
\parindent=10pt
\newdimen\nagykoz
\newdimen\kiskoz
\nagykoz=7pt
\kiskoz=2pt
\parskip=\nagykoz
\baselineskip=12.7pt


\loadeufm \loadmsam \loadmsbm

\font\vastag=cmssbx10
\font\drot=cmssdc10
\font\vekony=cmss10
\font\vekonydolt=cmssi10
\font\cimbetu=cmssbx10 scaled \magstep1
\font\szerzobetu=cmss10

\font\scVIII=cmcsc8
\font\rmVIII=cmr8
\font\itVIII=cmti8
\font\bfVIII=cmbx8
\font\ttVIII=cmtt8

\def\cim#1{{\centerline{\cimbetu#1}}}
\def\szerzo#1{{\vskip0.3truein\centerline{\szerzobetu#1}}}

\def\tetel#1#2{{{\drot#1}{\it\szukebb~#2\tagabb}}}
\long\def\biz#1#2{{{\vekony#1} #2}}
\def\kiemel#1{{\vekonydolt#1\/}}
\long\def\absztrakt#1#2{{\vskip0.4truein{\vekony#1} #2\vskip0.5truein}}
\def\szukebb{\parskip=\kiskoz}
\def\tagabb{\parskip=\nagykoz}
\def\vonal{{\vrule height 0.2pt depth 0.2pt width 0.5truein}}

\def\CC{{\Bbb C}}

\def\sbmfd{{Stein Banach manifold}}
\def\bmfd{{Banach manifold}}
\def\bsmfd{{Banach submanifold}}

\def\re{\hbox{\rm Re}}

\def\cts{{continuous}}
\def\bdd{{bounded}}

\def\fdml{{finite dimensional}}

\def\psh{{plurisubharmonic}}
\def\pshdom{{\psh\ domination}}
\def\holodom{{\holo\ domination}}

\def\st{{such that}}

\def\cln{{:}\;}
\def\iy{\infty}

\def\<{{\langle}}
\def\>{{\rangle}}

\def\RR {{\Bbb R}}

\def\OO {{\Cal O}}

\def\fn{func\-tion}
\def\fns{func\-tions}
\def\holo{hol\-o\-mor\-phic}

\def\mfd{manifold}
\def\smfd{submanifold}
\def\cpx{complex}
\def\cpt{compact}

\def\nbd{neighbor\-hood}

\def\bspc{Banach space}

\def\Da{{\Delta}}

\def\Th{{Theorem}}
\def\th{{theorem}}

\def\t#1{{\Th~#1}}

\newcount\minute    
\newcount\hour      
\newcount\hourMins  
%
%
\def\now%
{
%
  \minute=\time    
  \hour=\time \divide \hour by 60 
  \hourMins=\hour \multiply\hourMins by 60
  \advance\minute by -\hourMins 
  \zeroPadTwo{\the\hour}:\zeroPadTwo{\the\minute}%
}
%
%
\def\timestamp%
{
  \today\ \now
}
%
%
\def\today%
{
  \the\year-\zeroPadTwo{\the\month}-\zeroPadTwo{\the\day}%
}
%
%
\def\zeroPadTwo#1%
{
%
  \ifnum #1<10 0\fi    
  #1
}

{\phantom.}
\vskip0.5truein
\cim{ON COMPLEX BANACH MANIFOLDS}
\vskip0.2truein
\cim{SIMILAR TO STEIN MANIFOLDS}
\szerzo{Imre Patyi\plainfootnote{${}^*$}{\rmVIII 
Supported in part by NSF grant 
DMS 0600059.}}
\absztrakt{ABSTRACT.}{We give an abstract definition, similar to the
	axioms of a Stein manifold, of a class of complex Banach manifolds
	in such a way that a manifold belongs to the class if and only if it is
	biholomorphic to a closed split complex Banach submanifold of a separable
	Banach space.

	 MSC 2010: 46G20 (32Q28, 32Q35, 32Q40, 32T35, 32L20) 

         Key words: 
	 complex Banach manifold, 
	 holomorphic embedding,
	 sheaf cohomology.


}



\def\rD{D}
\def\rDPV{DPV}
\def\rF{F}
\def\rLA{L1}
\def\rLB{L}
\def\rLP{LP}
\def\rLf{Lf}
\def\rM{M}
\def\rP{P}
\def\rS{S}
\def\rZ{Z}


	 Stein manifolds can be characterized among \cpx\ \mfd{}s in various ways,
	including the two ways (I) and (II) below.
	 A paracompact second countable Hausdorff \cpx\ \mfd\ $M$ of pure dimension
	is a Stein manifold
	if and only if one and hence both of the following equivalent conditions (I)
	and (II) below hold.

	(I) (a) $M$ is \holo{}ally convex, i.e., if $K\subset M$ is \cpt,
	then its $\OO(M)$ \holo\ hull $\hat K$ is \cpt\ in $M$.
	(b) If $x\not=y$ in $M$, then there is an $f\in\OO(M)$ with $f(x)\not=f(y)$.
	(c) If $x\in M$, then there are an integer $n\ge1$ and a \holo\ \fn\
	$g\in\OO(M,\CC^n)$ that is
	a biholomorphism from an open \nbd\ $W$ of $x$ in $M$ to an open \nbd\
	$g(W)$ of $g(x)$ in $\CC^n$.

	(II) There is an $n\ge1$ \st\ $M$ is biholomorphic to a closed \cpx\ \smfd\
	$M'$ of $\CC^n$.

	 Let $X$ be a separable \bspc, and $M$ a paracompact second countable Hausdorff
	\cpx\ \bmfd\ modelled on $X$.
	 We call $M$ a \kiemel{\sbmfd\ modelled on $X$} if (i-iv) below hold.

	 (i) Holomorphic domination is possible in $M$, i.e., if $u\cln M\to\RR$
	is any locally upper bounded \fn, then there are a \bspc\ $Z$ and a \holo\
	\fn\ $h\cln M\to Z$ with $u(x)<\|h(x)\|$ for all $x\in M$.

	 (ii) There are open sets $U_n,V_n\subset M$, and \holo\ \fns\ $f_n\in\OO(M)$,
	$n\ge1$, \st\ $\bigcup_{n=1}^\iy(U_n\times V_n)=(M\times M)\setminus\Da_M$,
	where $\Da_M=\{(x,x)\cln x\in M\}$ is the diagonal of $M\times M$, and
	$f_n(U_n)$ and $f_n(V_n)$ are disjoint sets in $\CC$ for all $n\ge1$.

	 (iii) There are open sets $W_n\subset M$ and \holo\ maps $g_n\in\OO(M,X)$
	for $n\ge1$ \st\ $\bigcup_{n=1}^\iy W_n=M$ and $g_n|W_n$ is a biholomorphism from
	$W_n$ onto an open set $g_n(W_n)$ in $X$.

	 (iv) There are open sets $G_k\subset M$, $k\ge1$, with $\bigcup_{k=1}^\iy G_k=M$
	and $\sup_{x\in G_k}(|f_n(x)|+\|g_n(x)\|)<\iy$ for all $k,n\ge1$, where $(f_n)$
	and $(g_n)$ are as in (ii) and (iii).

	 If $M$ is \fdml, then it is easy to see that (i-iii) together are equivalent 
	to (I), and (iv) is vacuous, since if $G_k$, $k\ge1$, is an exhaustion of $M$
	by precompact open sets $G_k$, then any \cts\ \fn\ $|f_n(x)|+\|g_n(x)\|$ on $M$
	is \bdd\ on $G_k$ for $k,n\ge1$.
	 Thus if $M$ is \fdml, then (i-iv) together are equivalent to $M$ being a Stein
	\mfd.

\tetel{\t1.}{Let $X$ be a separable \bspc, and $M$ a paracompact second countable 
	Hausdorff \cpx\ \bmfd\ modelled on $X$.
	 Then $M$ is a \sbmfd\ modelled on $X$ if and only if there is a separable
	\bspc\ $X'$ \st\ $M$ is bi\holo\ to a closed split \cpx\ \bsmfd\ $M'$ of $X'$.
}

	 Here \bmfd{}s and \bsmfd{}s are understood in terms of bi\holo{}ally related
	charts, and a \bsmfd\ is called \kiemel{split} if each of its tangent spaces
	has a direct complement in the ambient \bspc.
	 Clearly, a \cpx\ \bsmfd\ $M$ of $X$ is split if and only if near each point
	$x_0\in M$ it is possible to split $X$ as a 
	direct sum $X=X'\times X''$ of closed linear subspaces $X',X''$ of $X$
	\st\ with $x_0=(x'_0,x''_0)$ and $x=(x',x'')$
	we can write $M$ as the graph $x''=m(x')$
	of a \holo\ \fn\ $m$ from an open \nbd\ of $x'_0$ in $X'$ to $X''$, where
	$x''_0=m(x'_0)$.

\biz{Proof.}{Suppose first that $M$ is bi\holo\ to an $M'$ and verify that
	$M$ satisfies (i-iv).
	 It is enough to show that $M'$ does.

	 As \holodom\ is possible in $X'$ by [\rP], and thus also in $M'$,
	since $M'$ is closed in $X'$, (i) is true.
	 We define some linear \fns\ $f_n\cln X'\to\CC$ and $g_n\cln X'\to X$ for
	$n\ge1$ whose restrictions to $M'$ will do the job.
	 For linear \fns\ (iv) is automatic: we can let $G_k$ be the intersection
	of $M'$ with the open ball $\|x\|<k$ in $X'$ and write
	$|f_n(x)|+\|g_n(x)\|\le(\|f_n\|+\|g_n\|)\|x\|\le(\|f_n\|+\|g_n\|)k<\iy$
	for $x\in G_k$ and $n,k\ge1$.

	 If $x\not=y$ in $M'$, then $x-y\not=0$ in $X'$ and the Hahn--Banach \th\ 
	gives us a \cpx\ linear functional $f_{xy}\in(X')^*$ of norm $1$ with
	$\re\,f_{xy}(x-y)=\|x-y\|>0$.
	 Let $U_{xy}=\{z\in X'\cln-\frac12\|x-y\|+\re\,f_{xy}(x)<\re\,f_{xy}(z)\}$
	and $V_{xy}=\{z\in X'\cln\re\,f_{xy}(z)<\frac12\|x-y\|+\re\,f_{xy}(y)\}$.
	 Then $x\in U_{xy}$, $y\in V_{xy}$, and their images $f_{xy}(U_{xy})$
	and $f_{xy}(V_{xy})$ are disjoint since they are the half planes
	$-\frac12\|x-y\|+\re\,f_{xy}(x)<\re\,w$, $\re\,w<\frac12\|x-y\|+\re\,f_{xy}(y)$,
	which are clearly disjoint since 
	$-\frac12\|x-y\|+\re\,f_{xy}(x)=\frac12\|x-y\|+\re\,f_{xy}(y)$.

	 Fix any point $x_0\in M'$ and denote its complex tangent space $T_{x_0}M'$
	by $X$ and regard it as a closed linear subspace of $X'$.
	 If $x\in M'$, then the complex tangent space $T_xM'$ and $X$ are linearly
	isomorphic via a bounded linear map $i_x\cln T_xM'\to X$, and there is a
	bounded linear projection $p_x\cln X'=T_xX'\to T_xM'$.
	 Thus the linear map $g_x\cln X'\to X$ given by $g_x(y)=i_x(p_x(y))$ for
	$y\in X'$ satisfies that $(dg_x)(x)y=i_x(y)$ for $y\in T_xM'$, i.e.,
	$(dg_x)(x)$ is a linear isomorphism from $T_xM'$ onto $X$.
	 By the inverse function \th\ $g_x$ is bi\holo\ from an open \nbd\ $W_x$ of
	$x$ in $M'$ to an open \nbd\ $g_x(W_x)$ of $g_x(x)=0$ in $X$.

	 By Lindel\"of's theorem in the second countable (separable metric) spaces 
	$(M'\times M')\setminus\Da_{M'}$ and $M'$ the open coverings
	$U_{xy}\times V_{xy}$, $(x,y)\in(M'\times M')\setminus\Da_{M'}$, and
	$W_x$, $x\in M'$, can be reduced to countable subcoverings $U_n\times V_n$,
	$W_n$, where $U_n=U_{x_ny_n}$, $V_n=V_{x_ny_n}$, and $W_n=W_{x'_n}$ for $n\ge1$.
	 Thus the \fns\ $f_n=f_{x_ny_n}$, $g_n=g_{x'_n}$, $n\ge1$, do the job.

	 Conversely, assume that $M$ satisfies (i-iv) and embed $M$ bi\holo{}ally
	as $M'$ into a separable \bspc\ $X'$.

	 If $i\ge1$, then let $C_i=L_i=1+\sup\{|f_n(x)|+\|g_n(x)\|\cln 1\le k,n\le i,
	x\in G_k\}$.
	 So if $k,n\ge1$, and $x\in G_k$, then $|f_n(x)|+\|g_n(x)\|\le C_kL_n$.
	 Thus upon replacing $f_n$ by $f_n/(L_n2^n)$ and $g_n$ by $g_n/(L_n2^n)$,
	we obtain new \fns\ again to be called $f_n$, $g_n$ that satisfy (ii), (iii),
	and the slightly strengthened version 
	$\sup_{x\in G_k}(|f_n(x)|+\|g_n(x)\|)<C_k/2^n$, $k,n\ge1$, of (iv).

	 The covering $W_n$, $n\ge1$, of the paracompact space $M$ has a locally
	finite refinement, which by Lindel\"of's \th\ can be taken to be countable,
	and can be shrunk since a paracompact Hausdorff space $M$ is normal.
	 There are open sets $M_n\subset M$, $n\ge1$, with $\bigcup_{n=1}^\iy M_n=M$,
	and for each $n\ge1$ there is an index $j(n)\ge1$ with the closure
	$\overline{M_n}\subset W_{j(n)}$.
	 Define $u\cln M\to\RR$ by $u(x)=\inf\{n\ge1\cln x\in M_n\}$.
	 Then $u$ is locally upper bounded on $M$ since $u\le n$ on the open set $M_n$.

	 By assumption (i) on \holodom\ there are a \bspc\ $Z$ and a \holo\ \fn\
	$h\in\OO(M,Z)$ with $u(x)<\|h(x)\|$ for $x\in M$.
	 As $Z'=\overline{\text{\rm span}}\{h(x)\cln x\in M\}$ is a separable \bspc,
	and as any separable \bspc\ can be embedded into $C[0,1]$ we can replace the
	\bspc\ $Z$ by the separable space $Z=C[0,1]$ endowed with the sup norm.

	 Define a \bspc\ $X'$ by $X'=Z\times\ell_1\times\ell_1(X)$, where $\ell_1$
	and $\ell_1(X)$ denote the spaces of summable sequences in $\CC$ and in $X$.
	 Let us write the variable $y$ in $X'$ as $y=(y',y''_n,y'''_n)$,
	where $y''=(y''_n)\in\ell_1$ and $y'''=(y'''_n)\in\ell_1(X)$.
	 Clearly, $X'$ is a separable \bspc, being the product of three such spaces.

	 Define the map $\Phi\cln M\to X'$ given by $y=\Phi(x)$, where
$$
	\cases
	y'\phantom{''_n}=h(x)\cr
	y''_n\phantom{'}=f_n(x)\cr
	y'''_n=g_n(x)\cr
	\endcases,
	\qquad n\ge1.
 $$
	 If $k\ge1$ and $x\in G_k$, then we have
	$\sum_{n=1}^\iy(|f_n(x)|+\|g_n(x)\|)\le\sum_{n=1}^\iy C_k/2^n\le C_k<\iy$.
	 Thus $\Phi$ is \holo.

	 Our $\Phi$ is injective, since if $x\not=y$ in $M$, then there is an index
	$n\ge1$ with $x\in U_n$ and $y\in V_n$, so $f_n(x)\not=f_n(y)$, and even more
	so $\Phi(x)\not=\Phi(y)$.

	 We claim that the set $M'=\Phi(M)$ is closed in $X'$.
	 Indeed, suppose that $y_i=\Phi(x_i)$, $x_i\in M$, converges 
	$y_i\to y$ in the norm in $X'$ as $i\to\iy$ to an
	element $y\in X'$, we must show that there is an $x\in M$ with $y=\Phi(x)$,
	i.e., $y\in M'$.
	 As $y'_i=h(x_i)$, $i\ge1$, is bounded, being convergent, 
	there is an index $b\ge1$,
	with $\|h(x_i)\|\le b$ for $i\ge1$, i.e., $u(x_i)<\|h(x_i)\|\le b$, or, 
	$x_i\in M_1\cup\ldots\cup M_b$ for $i\ge1$.
	 By the pigeon hole principle there are an index $a$ with $1\le a\le b$
	and an infinite set $I$ of indices $i$ \st\ $x_i\in M_a$ for all $i\in I$.
	 As $\overline{M_a}\subset W_{j(a)}$, and $g_{j(a)}$ is bi\holo\ on $W_{j(a)}$, 
	we see
	that $\Phi(x_i)$, $i\in I$, may converge only if the $x_i$, $i\in I$,
	converge in $\overline{M_a}$ to one of its elements $x\in\overline{M_a}\subset
	W_{j(a)}$.
	 Thus $y_i=\Phi(x_i)\to\Phi(x)=y$ as $i\to\iy$ in $I$.
	 As $M'$ contains $y$, it is closed.

	 If $y_0=\Phi(x_0)$ in $M'$, then there is an index $n\ge1$ with
	$x_0\in W_n$.
	 So $y'''_n=g_n(x)$ is bi\holo\ from a connected open set $W'_n$ with
	$x_0\in W'_n\subset\overline{W'_n}\subset W_n$ 
	to a connected open set $g_n(W'_n)$ in $X$.
	 Then the connected component of the set $M'\cap\{y'''_n\in g_n(W'_n)\}$
	that contains the point $y_0$ equals the graph
$$
	\cases
	y'''_n=y'''_n\cr
	y'\phantom{_n''}=h(g_n^{-1}(y'''_n))\cr
	y''_n\phantom{'}=f_n(g_n^{-1}(y'''_n))\cr
	y''_\nu\phantom{'}=f_\nu(g_n^{-1}(y'''_n))\cr
	y'''_\nu=g_\nu(g_n^{-1}(y'''_n))\cr
	\endcases,
	\qquad\nu\not=n,
 $$
 	of a \holo\ map $y'''_n\mapsto(y',y''_n,y'''_n,y''_\nu,y'''_\nu)$,
	$\nu\not=n$, from $W'_n$ to the \bspc\ $X'\cap\{y'''_n=0\}$.

	 Thus $M'$ is a closed split \cpx\ \bsmfd\ of $X'$ and $\Phi\cln M\to M'$ is
	a biholomorphism.
	 QED.
}

	 The most substantial part of the above proof is to show that \holodom\
	is possible on a separable \bspc.
	 That was done in the reference [\rP] based upon the work of Lempert
	in [\rLB].
	 It might be possible to weaken the axioms (i-iv) perhaps by dropping (iv)
	and replacing (i), that stands in for \holo\ convexity,
	by \pshdom, i.e., by requiring a \cts\ \psh\ \fn\ $\psi\cln M\to\RR$
	that dominates the given locally upper bounded \fn\ $u\cln M\to\RR$.
	 Nevertheless, axioms (i-iv) represent perhaps the ultimate axioms for
	``Stein Banach manifolds'' since any other system for which the desirable
	\t1 holds must be equivalent with (i-iv).
	 Most known methods of \pshdom\ also yield \holodom, and a `constructive'
	procedure for building the \fns\ $f_n,g_n$ in (ii) and (iii) is likely to
	produce \fns\ that also satisfy (iv).
	 The author doubts whether a successful ``Stein theory'' could be built
	up for nonseparable \bspc{}s and \bmfd{}s.
	 Even for separable \bmfd{}s it would be better to restrict attention to
	the ones modelled on separable \bspc{}s with the bounded approximation
	property (there are virtually no practical separable \bspc{}s that do
	not satisfy the bounded approximation property).
	 If $M$ is a \sbmfd\ modelled on such a \bspc, then the sheaf cohomology
	group $H^q(M,S)$ vanishes if $q\ge1$ and $S\to M$ is a so-called cohesive
	sheaf defined in [\rLP] by Lempert et al.
	 The question arises whether $M$ is a \sbmfd\ if $H^q(M,S)=0$ for all
	$q\ge1$ and all cohesive sheaves $S\to M$.
	 If $M$ is an open subset of a separable \bspc\ with the bounded approximation
	property, then the answer is Yes.

\vskip0.10truein
\centerline{\scVIII References}
\vskip0.10truein
\baselineskip=11pt
\parskip=3pt
\frenchspacing
{\rmVIII

\comment
	[\rD] Dineen, S.,
	{\itVIII 
	Complex analysis in infinite dimensional spaces},
	Springer-Verlag, 
	London,
	(1999).

	[\rDPV] \vonal, Patyi, I., Venkova, M.,
	{\itVIII 
	Inverses depending holomorphically on a parameter in a Banach space},
	J. Funct. Anal., 
	{\bfVIII 237} (2006), no. 1, 338--349.

	[\rF] Forstneri\v{c}, F.,
	{\itVIII
	Noncritical holomorphic functions on Stein manifolds},
	Acta Math.,
	{\bfVIII 191} (2003), 143--189.

	[\rLf] Laufer, H.\,B.,
	{\itVIII
	On the infinite dimensionality of the Dolbeault cohomology groups},
	Proc. Amer. Math. Soc.,
	{\bfVIII 52} No. 1 (1975), 293--296.

	[\rLA] Lempert,~L.,
	{\itVIII
	The Dolbeault complex in infinite dimensions~I},
	J. Amer. Math. Soc., {\bfVIII 11} (1998), 485--520.
\endcomment

	[\rLB] Lempert, L.,
	{\itVIII
	Plurisubharmonic domination},
	J. Amer. Math. Soc.,
	{\bfVIII 17}
	(2004),
	361--372.

	[\rLP] \vonal, Patyi,~I.,
	{\itVIII
	Analytic sheaves in Banach spaces},
	Ann. Sci. \'Ecole Norm. Sup., 
	S\'er. 4,
	{\bfVIII 40}
	(2007),
	453--486.

\comment
	[\rM] Mujica, J.,
	{\itVIII
	 Complex analysis in Banach spaces},
	North--Holland, Amsterdam,
	(1986).

        [\rP] Patyi, I.,
        {\itVIII
        On holomorphic Banach vector bundles over Banach spaces},
	Math. Ann.,
	{\bfVIII 341} (2008), no. 2, 455--482.

	[\rS] Schottenloher, M.,
	{\itVIII
	Spectrum and envelope of holomorphy for infinite-dimensional Riemann domains},
	Math. Ann.,
	{\bfVIII 263} (1983), no. 2, 213--219. 

	[\rZ] Zerhusen,~A.\,B.,
	{\itVIII
	Embeddings of pseudoconvex domains in certain Banach spaces},
	Math. Ann.,
	{\bfVIII 336}
	(2006), no. 2,
	269--280.
\endcomment

	[\rP] Patyi, I.,
	{\itVIII
	On holomorphic domination, I},
	arXiv:0910.0476.
}
\vskip0.20truein
\centerline{\vastag*~***~*}
\vskip0.15truein
{\scVIII
        Imre Patyi,
        Department of Mathematics and Statistics,
        Georgia State University,
        Atlanta, GA 30303-3083, USA,
        {\ttVIII ipatyi\@gsu.edu}
}
\bye